\documentclass[12pt,reqno]{amsart}
\usepackage{amssymb,verbatim,enumerate,ifthen,graphics,graphicx}
\usepackage[mathscr]{eucal}
\usepackage[utf8]{inputenc}
\usepackage[T1]{fontenc}
\newcommand{\R}{\mathbb{R}}
\newcommand{\N}{\mathbb{N}}
\newcommand{\Z}{\mathbb{Z}}
\renewcommand{\P}{\mathscr{P}}
\newcommand{\U}{\mathscr{U}}

\newcommand{\dist}{\mathop{\hbox{\rm dist}}\nolimits}

\newcommand{\rec}{\mathop{\hbox{\rm rec}}\nolimits}

\newcommand{\cl}{\mathop{\hbox{\rm cl}}\nolimits}

\newtheorem{theorem}{Theorem}

\newtheorem*{theorem*}{Theorem}
\def\Thm#1#2{\ifthenelse{\equal{#1}{*}}{\begin{theorem*}#2\end{theorem*}}
  {\begin{theorem}\label{T#1}#2\end{theorem}}}
\newtheorem{Atheorem}{Theorem}

\def\thm#1{Theorem~\ref{T#1}}
\newtheorem{proposition}[theorem]{Proposition}
\newtheorem*{proposition*}{Proposition}
\def\Prp#1#2{\ifthenelse{\equal{#1}{*}}{\begin{proposition*}#2\end{proposition*}}
             {\begin{proposition}\label{P#1}#2\end{proposition}}}
\def\prp#1{Proposition~\ref{P#1}}
\newtheorem{corollary}[theorem]{Corollary}
\newtheorem*{corollary*}{Corollary}
\def\Cor#1#2{\ifthenelse{\equal{#1}{*}}{\begin{corollary*}#2\end{corollary*}}
             {\begin{corollary}\label{C#1}#2\end{corollary}}}
\def\cor#1{Corollary~\ref{C#1}}
\newtheorem{lemma}[theorem]{Lemma}
\newtheorem*{lemma*}{Lemma}
\def\Lem#1#2{\ifthenelse{\equal{#1}{*}}{\begin{lemma*}#2\end{lemma*}}
             {\begin{lemma}\label{L#1}#2\end{lemma}}}
\def\lem#1{Lemma~\ref{L#1}}
\newtheorem{example}[theorem]{Example}
\newtheorem*{example*}{Example}
\long\def\Exa#1#2{\ifthenelse{\equal{#1}{*}}{\begin{example*}\rm #2\end{example*}}
            {\begin{example}\label{Ex#1}\rm #2\end{example}}}

\newtheorem{problem}[subsection]{Problem}

\theoremstyle{definition}
\newtheorem{definition}[theorem]{Definition}
\newtheorem*{definition*}{Definition}
\def\Defi#1#2{\ifthenelse{\equal{#1}{*}}{\begin{definition*}#2\end{definition*}}
      {\begin{definition}\label{D#1}#2\end{definition}}}

\newtheorem{remark}[theorem]{Remark}
\newtheorem*{remark*}{Remark}
\long\def\Rem#1#2{\ifthenelse{\equal{#1}{*}}{\begin{remark*}#2\end{remark*}}
             {\begin{remark}\label{R#1}#2\end{remark}}}
\def\rem#1{Remark~\ref{R#1}}
\def\eq#1{{\rm(\ref{E#1})}}
\def\Eq#1#2{\ifthenelse{\equal{#1}{*}}
  {\begin{equation*}\begin{aligned}[]#2\end{aligned}\end{equation*}}
  {\begin{equation}\begin{aligned}[]\label{E#1}#2\end{aligned}\end{equation}}}

\allowdisplaybreaks

\begin{document}

\begin{flushright}
\textit{To appear in Publ. Math. Debrecen} 
\end{flushright}
\vspace{5mm}

\date{\today}

\title[On strongly and approximately convex and concave maps]
{Bernstein--Doetsch type theorems for set-valued maps of strongly and approximately convex and concave type}

\author[C. Gonz\'alez]{Carlos Gonz\'alez}
\address{Escuela de Matematicas, Universidad Central de Venezuela, Caracas, Venezuela}
\email{carlosl.gonzalez@ciens.ucv.ve}

\author[K. Nikodem]{Kazimierz Nikodem}
\address{Department of Mathematics and Computer Science, University of Bielsko-Bia{\l}a, ul.\
Willowa~2, 43-309 Bielsko-Bia{\l}a, Poland}
\email{knikodem@ath.bielsko.pl}

\author[Zs. P\'ales]{Zsolt P\'ales}
\address{Institute of Mathematics, University of Debrecen, H-4010 Debrecen, Pf.\ 12, Hungary}
\email{pales@science.unideb.hu}

\author[G. Roa]{Gari Roa}
\address{Escuela de Matematicas, Universidad Central de Venezuela, Caracas, Venezuela}
\email{gariroa@gmail.com}

\dedicatory{This paper is dedicated to the 90th birthday of Professor Lajos Tamássy.}


\thanks{This research of the third author was realized in the frames of T\'AMOP 4.2.4. A/2-11-1-2012-0001 
”National Excellence Program – Elaborating and operating an inland student and researcher personal
support system”. The project was subsidized by the European Union and co-financed by the European 
Social Fund. This research of the third author was also supported by the Hungarian Scientific 
Research Fund (OTKA) Grant NK 81402.}

\begin{abstract}
In this paper, we investigate properties of set-valued mappings that establish connection between
the values of this map at two arbitrary points of the domain and the value at their midpoint. 
Such properties are, for instance, Jensen convexity/concavity, $K$-Jensen convexity/concavity 
(where $K$ is the set of nonnegative elements of an ordered vector space), and approximate/strong 
$K$-Jensen convexity/concavity. Assuming weak but natural regularity assumptions on the set-valued map, 
our main purpose is to deduce the convexity/concavity consequences of these properties in the appropriate 
sense. Our two main theorems will generalize most of the known results in this field, in particular the 
celebrated Bernstein--Doetsch Theorem from 1915, and thus they offer a unified view of these theories. 
\end{abstract}

\maketitle

\section{Introduction}

The classical theorem of Bernstein and Doetsch \cite{BerDoe15} published almost one hundred years ago plays, undoubtedly, a fundamental role in the theory of convexity. This theorem asserts that if $f:D\to\R$ is a 
Jensen convex function (where $D$ is a real interval, or more generally, $D$ is a convex subset of a linear space $X$), i.e.,
\Eq{JC}{
  f\Big(\frac{x + y}{2}\Big)\leq \frac{f(x) + f(y)}{2} \qquad(x,y\in D)
}
and it is also locally upper bounded, then it must be convex over $D$. If $(-f)$ is Jensen convex, then $f$ is frequently called Jensen concave, the results for these functions are analogous under the assumption of local lower boundedness. The above theorem was seminal and was applied and generalized in many useful ways to various important circumstances that we briefly describe in what follows.  

When the co-domain of $f$ is an ordered vector space $Y$, i.e., the set $K$ of the nonnegative elements in $Y$ 
forms a convex cone, then one can define $K$-Jensen convexity of $f:D\to Y$ by
\Eq{JCK}{
  \frac{f(x) + f(y)}{2}\in f\Big(\frac{x + y}{2}\Big)+K \qquad(x,y\in D).
}
In particular, if $Y=\R$ and $K=\R_+$, then \eq{JCK} is equivalent to \eq{JC}. Analogously, one can introduce
$K$-convexity. Extensions of the Bernstein--Doetsch Theorem to this setting were formulated by Trudzik \cite{Tru84}. Functions $f:D\to Y$ satisfying the inclusion
\Eq{JCVK}{
  f\Big(\frac{x + y}{2}\Big)\in \frac{f(x) + f(y)}{2} + K \qquad(x,y\in D)
}
are called $K$-Jensen concave. Obviously, this holds if and only if $(-f)$ is $K$-Jensen convex (or if
$f$ is $(-K)$-Jensen convex). Thus, the results related to $K$-Jensen concavity of functions can always 
be derived from the statements established for $K$-Jensen convexity. (This, however, will not be the case for the 
set-valued setting.)

A further generalization step was to consider set-valued functions instead of single-valued ones.
A set-valued map $F:D\to 2^Y$ is called $K$-Jensen convex if the inclusion
\Eq{SVJCK}{
  \frac{F(x) + F(y)}{2}\subseteq F\Big(\frac{x + y}{2}\Big)+K \qquad(x,y\in D)
}
holds. Observe that if $F$ is of the form $F(x)=\{f(x)\}$ for some function $f:D\to Y$, then \eq{SVJCK}
is equivalent to \eq{JCK}, and hence \eq{SVJCK} generalizes \eq{JCK}. Bernstein--Doetsch type 
results for this setting have been obtained by 
Averna, Cardinali, Nikodem, and Papalini \cite{AveCar90,CarNikPap93,Nik86,Nik87a,Nik87c,Nik89,Pap90} 
and by Borwein \cite{Bor77}. The notion of the $K$-Jensen concavity of a set-valued function
$F:D\to 2^Y$, analogously to \eq{JCVK}, is defined by 
\Eq{SVJCVK}{
  F\Big(\frac{x + y}{2}\Big)\subseteq \frac{F(x) + F(y)}{2} + K \qquad(x,y\in D).
}
Observe that, in general, the $K$-Jensen concavity of $F$ is \textit{not} equivalent to the $K$-Jensen 
convexity of $(-F)$. Thus, starting from the set-valued setting, the cases of convexity and concavity 
need separate (but similar) considerations.

Another chain of generalizations of the Bernstein--Doetsch Theorem emerged from the paper of Ng and Nikodem \cite{NgNik93} in the context of approximate convexity. They proved if $f:D\to\R$ is $\varepsilon$-Jensen convex for some $\varepsilon\geq0$, i.e.,
\Eq{JCe}{
  f\Big(\frac{x + y}{2}\Big)\leq \frac{f(x) + f(y)}{2} + \varepsilon \qquad(x,y\in D),
}
and if $f$ is also locally upper bounded, then it is $2\varepsilon$-convex, i.e.,
\Eq{Ce}{
  f(tx+(1-t)y)\leq tf(x)+(1-t)f(y)+2\varepsilon \qquad(x,y\in D,\,t\in[0,1]).
}
This result was also independently established by Laczkovich \cite{Lac99}.
The set-valued and $K$-Jensen convex/concave variant of this theorem will be formulated as a consequence of 
our results in the sequel. 

Considering more general error terms than constant ones, Házy and Páles 
\cite{HazPal04} investigated the following approximate Jensen type inequality:
\Eq{JCHP}{
  f\Big(\frac{x + y}{2}\Big)\leq \frac{f(x) + f(y)}{2} + \varepsilon\|x-y\| \qquad(x,y\in D),
}
provided that $D$ is a subset of a normed space $X$ and $f$ is real valued. They showed, under the usual local upper-boundedness condition, that \eq{JCHP} implies
\Eq{CHP}{
  f(tx+(1-t)y)\leq tf(x)+(1-t)f(y)+2\varepsilon T(t)\|x-y\| \qquad(x,y\in D,\,t\in[0,1]),
}
where the function $T:\R\to\R$, the so-called Takagi function, is defined by
\Eq{Tak}{
  T(t):=\sum_{n=0}^\infty \frac1{2^n}\dist(2^nt,\Z).
}
Results, extending this approach to more general error terms and also to convexity concepts related to Chebyshev systems, have recently been obtained by Házy, Makó and Páles \cite{Haz05a,Haz07b,HazPal05,HazPal09,MakPal10b,MakPal12b,MakPal12c,MakPal13b} and by
Mureńko, Ja. Tabor, Jó. Tabor, and Żoldak \cite{MurTabTab12,TabTab09b,TabTab09a,TabTabZol10b,TabTabZol10a}.
For the set-valued and $K$-Jensen convex/concave setting more general statements will be formulated as direct 
consequences of our two main results below. 

Finally, we mention the notion of strong convexity, which in some sense, is opposite to approximate convexity.
Following Polyak \cite{Pol66}, a function $f:D\to\R$ is called strongly Jensen convex with modulus $\varepsilon\geq0$ if
\Eq{SJC}{
  f\Big(\frac{x + y}{2}\Big)\leq \frac{f(x) + f(y)}{2} - \frac{\varepsilon}{4}\|x-y\|^2 \qquad(x,y\in D).
}
Assuming local upper boundedness of $f$, Azócar, Gimenez, Nikodem and Sanchez in \cite{AzoGimNikSan11} showed that the above strong Jensen-convexity property implies that $f$ is strongly convex with modulus $\varepsilon$, i.e., 
\Eq{SC}{
  f(tx+(1-t)y)\leq tf(x)+(1-t)f(y) - \varepsilon t(1-t)\|x-y\|^2 \qquad(x,y\in D,\,t\in[0,1]).
}
The set-valued Jensen convex variant of this result was established by Leiva, Merentes, Nikodem, and Sanchez \cite{LeiMerNikSan13}. More abstract and powerful results will be the corollaries of our two main theorems.

\section{Terminology and Auxiliary Results}

Let $X$ be a Hausdorff topological linear space. The family of all open sets containing the origin of $X$
will be denoted by $\U=\U(X)$. The closure of a subset $H\subseteq X$ will be denoted either by 
$\overline{H}$ or by $\cl(H)$. We will frequently use the fact that $\overline{H}=\bigcap_{U\in\U(X)}(U+H)$.
An immediate consequence of this formula is that, for every pair of subsets $A,B\subseteq X$, we have
\Eq{cl+}{
 \overline{A}+\overline{B}\subseteq \overline{A+B}=\overline{\overline{A}+\overline{B}}.
}
(See, e.g., \cite{Rud91}, \cite{Rad60}).
Indeed, to prove the inclusion in \eq{cl+}, let $U\in\U(X)$ and choose $V\in\U(X)$ such that 
$V+V\subseteq U$. Then $\overline{A}+\overline{B}\subseteq (V+A)+(V+B)\subseteq U+A+B$. 
Hence $\overline{A}+\overline{B}\subseteq \bigcap_{U\in\U(X)}(U+A+B)=\overline{A+B}$.
The equality in \eq{cl+} trivially follows from this inclusion and the monotonicity of the closure
operation. If either $\overline{A}$ or $\overline{B}$ is compact then, the set $\overline{A}+\overline{B}$
is closed and the inclusion in \eq{cl+} can be replaced by equality. 

A subset $H\subseteq X$ is called \textit{bounded} if, for any $U\in\U$, there exists $t>0$ such that
$H\subseteq tU$. It is well-known that the union of finitely many bounded sets is also bounded. 
The family of bounded sets is closed under algebraic addition and multiplication by scalars. 

Given two points $x,y\in X$, their convex hull, i.e., the \textit{segment} connecting them, 
will be denoted by $[x,y]$. \textit{Convexity of a set $H\subseteq X$} is equivalent to the 
property that $[x,y]\subseteq H$ whenever $x,y\in H$. A set $H\subseteq X$ is said to be 
\textit{closedly convex} if, for all $x,y\in H$, $[x,y]\subseteq\overline{H}$ holds. In fact, 
it is easy to see that $H$ is closedly convex if and only if $\overline{H}$ is convex.
A convex set which is closed with respect to multiplication by positive scalars is called a
\textit{convex cone}. A set $H\subseteq X$ is called \textit{star-shaped with respect to 
a point $p\in H$}, if for all $x\in H$, the segment $[x,p]$ is contained in $H$.

Given a convex cone $K\subseteq X$, we can define the \textit{ordering relation} $\leq_K$ as follows: 
\textit{$x\leq_K y$ holds if and only if $y-x\in K$}. This relation is transitive, and, provided
that $K\cap(-K)=\{0\}$, it is also reflexive and antisymmetric and then $\leq_K$ is a partial order on
$X$. 

Given a convex cone $K$, a subset $S\subseteq X$ is called \textit{$K$-lower bounded} if there exists a
bounded set $H\subseteq X$ such that $S\subseteq H+K$. If $S\subseteq \cl(H+K)$ for some bounded set $H$,
then we say that $S$ is \textit{closedly $K$-lower bounded}. If the space $X$ is locally bounded,
then these notions are equivalent to each other. One can easily see that the union of finitely many
(closedly) $K$-lower bounded sets is also (closedly) $K$-lower bounded. Furthermore, the family of
(closedly) $K$-lower bounded sets is closed under algebraic addition and multiplication by positive scalars.

A set $H\subseteq X$ is called \textit{$K$-convex} (resp.\ closedly \textit{$K$-convex}) if $[x,y]$ is 
contained in $H+K$ (resp.\ in $\cl(H+K)$) whenever $x,y\in H$. In general, using the inclusion in \eq{cl+}, 
one can see that the $\overline{K}$-convexity of $\overline{H}$ implies that $H$ is closedly $K$-convex, 
but the reversed implication may not be valid for all convex cones $K$.

Given a nonempty subset $H\subseteq X$, we can naturally attach a cone, called the
\textit{recession cone} of $H$ in the following manner:
\Eq{*}{
  \rec(H):=\{x\in X\mid tx+H\subseteq H\mbox{ for all }t\geq0\}.
}

It is elementary to check that $\rec(H)$ is a convex cone containing $0$.
The additional basic properties of the recession cone are summarized in the following lemma.

\Lem{rec}{Let $H\subseteq X$ be a nonempty set. Then
\begin{enumerate}[(i)]
 \item $\rec(H)$ is a convex cone containing $0$;
 \item $K=\rec(H)$ is the largest cone $K$ such that $K+H\subseteq H$ is valid;
 \item $\overline{\rec}(H)\subseteq\rec(\overline{H})$;
 \item for all $x\in X$, $t>0$, $\rec(x+tH)=\rec(H)$;
 \item for all nonempty sets $H_1,H_2\subseteq X$, $\rec(H_1)+\rec(H_2)\subseteq \rec(H_1+H_2)$.
\end{enumerate}}

\begin{proof}
The proofs of (i), (ii) and (iv) are immediate. 

Taking closure of both sides of the inclusion $\rec(H)+H\subseteq H$ and then using \eq{cl+}, it follows that $\overline{\rec}(H)+\overline{H}\subseteq\overline{H}$. Hence, using (ii), we conclude that $\overline{\rec}(H)\subseteq\rec(\overline{H})$ proving (iii).

Adding up the relations $\rec(H_1)+H_1\subseteq H_1$ and $\rec(H_2)+H_2\subseteq H_2$ side by side,
we get $(\rec(H_1)+\rec(H_2))+H_1+H_2\subseteq H_1+H_2$. Hence, using (ii), the inclusion stated in (v) follows.
\end{proof}

In order to prove our main results in the next section, we will also need the 
following lemma, which allows one to perform the limit related to inclusions.

\Lem{Inc}{Let $(A_n),(B_n)$ be nondecreasing sequences of nonempty subsets of $X$, let
$H\subseteq X$ be a nonempty bounded set, let $K\subseteq \overline{\rec}(B_0)$ be also 
nonempty, and let $\varepsilon_n>0$ be a null-sequence of real numbers. Assume that, for all
$n\geq0$,
\Eq{Incn}{
  A_n\subseteq \cl\big(\varepsilon_n H + K + B_n\big).
}
Then
\Eq{Inc}{
  \cl\bigg(\bigcup_{n=0}^{\infty}A_n\bigg)
  \subseteq \cl\bigg(\bigcup_{n=0}^{\infty}B_n\bigg).
}}

\begin{proof}
Let $U\subseteq X$ be any neighborhood of the origin $0\in X$ and let $V\subseteq X$ be an open
balanced neighborhood of $0$ such that $V+V+V\subseteq U$. The set $K$ being a subset of the closure
of the recession cone $\rec(B_n)$, we have that $K\subseteq V+\rec(B_n)$. By definition,
we also have that $\rec(B_n)+B_n\subseteq B_n$. Thus, for all $n\geq0$, it follows that
$K+B_n\subseteq V+B_n$.

On the other hand, since $H$ is bounded and $\varepsilon_n$ is a null-sequence, there exists an
integer $N$ such that, for $n\geq N$, $\varepsilon_n H\subseteq V.$ Therefore, applying \eq{Incn},
for $n\geq N$, we get that
\Eq{*}{
 A_n\subseteq \cl\big(\varepsilon_n H + K + B_n\big)
    \subseteq \cl\big(V+ V + B_n\big) 
    &\subseteq V+ V + V + B_n
    \subseteq  U + B_n
    \subseteq  U + \bigcup_{k=0}^\infty B_k.
}
Therefore,
\Eq{Inc2}{
 \bigcup_{n=N}^\infty A_n
 \subseteq U + \bigcup_{k=0}^\infty B_k. 
}
In view of the nondecreasingness of the sequence $(A_n)$, inclusion \eq{Inc2} is equivalent to
\Eq{*}{
 \bigcup_{n=0}^\infty A_n
 \subseteq U + \bigcup_{k=0}^\infty B_k. 
}
The above relation is valid for all $U\in\U(X)$, hence
\Eq{*}{
 \bigcup_{n=0}^\infty A_n
 \subseteq \bigcap_{U\in\U(X)}\bigg(U+\bigcup_{k=0}^\infty B_k\bigg)
 = \cl\bigg(\bigcup_{k=0}^\infty B_k\bigg). 
}
The right hand side of the above inclusion is closed, therefore \eq{Inc} follows immediately.
\end{proof}

Given another Hausdorff topological linear space $Y$, we will now introduce further notions
for set-valued maps from a subset of $X$ to $\P_0(Y)$, where $\P_0(Y)$ stands for the family of
nonempty subsets of the space $Y$. 

Let $D\subseteq X$ be a nonempty set and $K$ be a convex cone in $Y$. A set-valued function
$S:D\to\P_0(Y)$ is called \textit{locally closedly $K$-lower bounded} if, for any $x\in D$, there exist an
open set $U$ containing $x$ and a bounded set $H\subseteq Y$ such that $S(u)\subseteq \cl(H+K)$
holds for all $u\in U\cap D$. A set-valued function $S:D\to\P_0(Y)$ is called \textit{locally closedly 
weakly $K$-upper bounded} if, for any $x\in D$, there exist an open set $U$ containing $x$
and a bounded set $H\subseteq Y$ such that $0\in \cl(S(u)+H+K)$ for all $u\in U\cap D$.

\Lem{lbd}{Assume that $S:D\to\P_0(Y)$ is locally closedly $K$-lower bounded set-valued map. Then, for each
compact subset $C\subseteq D$, there exists a bounded set $H\subseteq Y$ such that, for all $x\in
C$, $S(x)\subseteq \cl(H+K)$.}

\begin{proof} Let $C\subseteq D$ be a nonempty compact set. Since $S$ is locally closedly $K$-lower
bounded, therefore, for each
$x\in C$, there exist an open set $U_x$ containing $x$ and a bounded set $H_x\subseteq Y$ such that
$S(u)\subseteq \cl(H_x+K)$ for all $u\in U_x\cap D$. The family of sets $\{U_x\mid x\in C\}$ is an open
covering for $C$, hence, by the compactness of $C$, there exists a finite subcovering
$\{U_{x_1},\dots,U_{x_n}\}$ of $C$. Let $H:=H_{x_1}\cup\cdots\cup H_{x_n}$. Then, $H$ is bounded
and, for all $x\in C$, there exists $i\in\{1,\dots,n\}$ such that $x\in U_{x_i}$. Hence
$S(x)\subseteq \cl(H_{x_i}+K)\subseteq \cl(H+K)$, which proves the assertion.
\end{proof}

\Lem{wubd}{Assume that $S:D\to\P_0(Y)$ is locally closedly weakly $K$-upper bounded set-valued map. 
Then, for each compact subset $C\subseteq D$, there exists a bounded set $H\subseteq Y$ such that, for all
$x\in C$, $0\in \cl(S(x)+H+K)$.}

\begin{proof} Let $C\subseteq D$ be a nonempty compact set. The set-valued map $S$ being locally 
closedly weakly $K$-upper bounded, it follows that, for each
$x\in C$, there exist an open set $U_x$ containing $x$ and a bounded set $H_x\subseteq Y$ such that
$0\in \cl(S(u)+ H_x+K)$ for all $u\in U_x$. The family of sets $\{U_x\mid x\in C\}$ is an open
covering for $C$, hence, by the compactness of $C$, there exists a finite subcovering
$\{U_{x_1},\dots,U_{x_n}\}$ of $C$. Let $H:=H_{x_1}\cup\cdots\cup H_{x_n}$. Then, $H$ is bounded
and, for all $x\in C$, there exists $i\in\{1,\dots,n\}$ such that $x\in U_{x_i}$. Hence $0\in
\cl(S(x)+H_{x_i}+K)\subseteq \cl(S(x)+H+K)$, which proves the assertion.
\end{proof}

\section{Takagi transformation of set-valued maps}

Assume now that $D\subseteq X$ is a star-shaped set and consider now a set-valued map
$S:D\to\P_0(Y)$ with the additional property that $0\in S(x)$ for all $x\in D$. For such a map,
we define $S^T:\R\times D\to Y$ by the following expression:
\Eq{Tak0}{
  S^T(t,x):=\cl\bigg(\bigcup_{n=0}^{\infty} \sum_{k=0}^{n} 
                 \frac{1}{2^k}S\big(2d_{\Z}(2^kt)x\big)\bigg)\qquad(t\in\R,\,x\in D),
}
where $d_{\Z}:\R\to\R$ is defined by
\Eq{*}{
  d_{\Z}(t):=\dist(t,\Z):=\inf_{k\in\Z}|k-t| \qquad(t\in\R).
}
The set-valued map $S^T$ will be called the \textit{Takagi transformation} of the map $S$ in the sequel.
The recession cone of a set-valued map $S$ is set by
\Eq{*}{
 \rec(S):=\bigcap_{x\in D} \rec S(x).
}

In the following lemma we establish the relationship between a set-valued map and its 
Takagi transformation.

\Lem{TT}{Let $D\subseteq X$ be a star-shaped set and $S:D\to\P_0(Y)$ be a set-valued map with 
the additional property that $0\in S(x)$ for all $x\in D$. Then
\Eq{TT1}{
   \cl(S(x))\subseteq S^T\big(\tfrac12,x\big) \qquad (x\in D).
}
If, in addition $S(0)\subseteq \overline\rec(S)$, then
\Eq{TT2}{
   \cl(S(x))= S^T\big(\tfrac12,x\big) \qquad (x\in D).
}}

\begin{proof} Observe that $d_{\Z}\big(\tfrac12\big)=\tfrac12$ and $d_{\Z}\big(2^k\cdot\tfrac12\big)=0$ for $k\in\N$.
Thus,
\Eq{*}{
  S^T\big(\tfrac12,x\big)
   =\cl\bigg(\bigcup_{n=0}^{\infty} \sum_{k=0}^{n}\frac{1}{2^k}S\big(2d_{\Z}(2^k\cdot\tfrac12)x\big)\bigg)
   =\cl\bigg(S(x)+\bigcup_{n=0}^{\infty} \sum_{k=1}^{n}\frac{1}{2^k}S(0)\bigg).
}
The property $0\in S(0)$, directly implies that \eq{TT1} holds. To prove \eq{TT2}, assume 
that $S(0)\subseteq \overline\rec(S)$. Then $S(0)\subseteq \overline\rec(S(x))\subseteq\rec(\overline{S(x)})$.
Since $\rec(\overline{S(x)})$ is a convex cone, we have that this set is closed under addition and 
multiplication by scalars. Thus, for all $n\in\N$,
\Eq{*}{
  \sum_{k=1}^{n}\frac{1}{2^k}S(0)
  \subseteq\sum_{k=1}^{n}\frac{1}{2^k}\rec(\overline{S(x)})
  \subseteq\rec(\overline{S(x)}).
}
Consequently, 
\Eq{*}{
\bigcup_{n=0}^{\infty} \sum_{k=1}^{n}\frac{1}{2^k}S(0)\subseteq\rec(\overline{S(x)}).
}
Therefore,
\Eq{*}{
  S^T\big(\tfrac12,x\big)
   =\cl\bigg(S(x)+\bigcup_{n=0}^{\infty}\sum_{k=1}^{n}\frac{1}{2^k}S(0)\bigg)
   \subseteq\cl\bigg(\overline{S(x)}+\rec(\overline{S(x)})\bigg)
   \subseteq \cl(S(x)),
}
which completes the proof of \eq{TT2}. 
\end{proof}

The following lemma implies that the Takagi transformation of a set-valued map which is constructed
as the product of an upper semicontinuous nonnegative scalar function and a convex subset of $Y$ is
the product of the Takagi transformation of the scalar function and the same set.

\Prp{Tak}{Let $D\subseteq X$ be a star-shaped set and $S_0\subseteq Y$ be a convex set 
containing $0\in Y$ and $K\subseteq Y$ be a convex cone. Let $\varphi:D\to \R_+$ be a locally upper
bounded nonnegative function. Define $S:D\to\P_0(Y)$ by $S(x):=K+\varphi(x)S_0$. Then
\Eq{Tak1}{
  S^T(t,x)=\cl\big(K+\varphi^T(t,x) S_0\big)
  \qquad(t\in\R,\,x\in D),
}
where
\Eq{Tak2}{  
  \varphi^T(t,x)=\sum_{n=0}^{\infty}\frac{1}{2^n}\varphi\big(2d_{\Z}(2^nt)x\big)
  \qquad(t\in\R,\,x\in D).
}
If, in addition, $\varphi(0)=0$, then
\Eq{Tak+}{
  \varphi^T\big(\tfrac12,x\big)=\varphi(x) \qquad\mbox{and}\qquad 
  S^T\big(\tfrac12,x\big)=\cl(K+\varphi(x)S_0)=\cl(S(x))\qquad(x\in D).
}}

\begin{proof} For $t\in\R$ and $n\geq 0$, we have that $0\leq 2d_{\Z}(2^nt)\leq 1$,
therefore $2d_{\Z}(2^nt)x\in[0,x]$. The function $\varphi$ being locally upper bounded on $D$, 
it is bounded from above over $[0,x]$ by some constant $M(x)$. Then we have that
\Eq{*}{
  \varphi^T(t,x)=\sum_{n=0}^{\infty}\frac{1}{2^n}\varphi\big(2d_{\Z}(2^nt)x\big) 
    \leq \sum_{n=0}^{\infty}\frac{1}{2^n}M(x)=2M(x) \qquad(t\in\R).
}

To prove \eq{Tak1}, fix $(t,x)\in\R\times D$. 

For the proof of the inclusion $\subseteq$ in \eq{Tak1}, we first show that
\Eq{Tak3}{
\bigcup_{n=0}^{\infty} \sum_{k=0}^{n}\frac{1}{2^k}S\big(2d_{\Z}(2^kt)x\big)
\subseteq K+\varphi^T(t,x) S_0.
}
Choose $y$ from the left hand side of inclusion \eq{Tak3} arbitrarily.
Then, using the definition of $S$ and the convexity of $S_0$, for some $n\geq 0$ we have that
\Eq{*}{
  y &\in\sum_{k=0}^{n} \frac{1}{2^k}S\big(2d_{\Z}(2^kt)x\big) 
      = \sum_{k=0}^{n} \bigg(K + \frac{1}{2^k}\varphi\big(2d_{\Z}(2^kt)x\big)S_0\bigg) \\
    & = K + \bigg(\sum_{k=0}^{n} \frac{1}{2^k}\varphi\big(2d_{\Z}(2^kt)x\big)\bigg)S_0
        +\bigg(\sum_{k=n+1}^{\infty}\frac{1}{2^k}\varphi\big(2d_{\Z}(2^kt)x\big)\bigg)\{0\} \\
    &\subseteq K + \bigg(\sum_{k=0}^{\infty} \frac{1}{2^k}\varphi\big(2d_{\Z}(2^kt)x\big)\bigg)S_0
      = K+\varphi^T(t,x) S_0.
}
Thus, inclusion \eq{Tak3} is proved. Taking the closures of both sided, 
the inclusion $\subseteq$ in \eq{Tak1} follows immediately. 

For the proof of the inclusion $\supseteq$ in \eq{Tak1}, by the closedness of $S^T(t,x)$,
it suffices to show that
\Eq{*}{
  K+\varphi^T(t,x)S_0\subseteq S^T(t,x).
}
Let $y\in K+\varphi^T(t,x)S_0$.
Then $y=u+\varphi^T(t,x)v$ for some elements $u\in K$, $v\in S_0$. Define the sequence
$(y_n)$ by 
\Eq{*}{
y_n:= u+\sum_{k=0}^{n} \frac{1}{2^k}\varphi\big(2d_{\Z}(2^kt)x\big)v.
}
Obviously, $y_n\to y$ as $n\to\infty$. On the other hand, for all $n\geq 0$
\Eq{*}{
 y_n\in K + \sum_{k=0}^{n} \frac{1}{2^k}\varphi\big(2d_{\Z}(2^kt)x\big)S_0
     &= \sum_{k=0}^{n} \frac{1}{2^k} \big(K+\varphi\big(2d_{\Z}(2^kt)x\big)S_0\big) \\
     &= \sum_{k=0}^{n} \frac{1}{2^k}S\big(2d_{\Z}(2^kt)x\big)\subseteq S^T(t,x).
}
By the closedness of $S^T(t,x)$, it follows that the limit $y$ of the sequence $(y_n)$ is
also contained in $S^T(t,x)$, which completes the proof of \eq{Tak1}.

In the case of $\varphi(0)=0$, the first equality in \eq{Tak+} is immediate, the second
equality is then a consequence of \eq{Tak1}. 
\end{proof}

\Cor{Tak}{Let $X$ be a normed space, $D\subseteq X$ be a star-shaped set, $S_0\subseteq Y$ be a convex set 
containing $0\in Y$, $K\subseteq Y$ be a convex cone, and $\alpha>0$. Define $S:D\to\P_0(Y)$ by
$S(x):=K+\|x\|^\alpha S_0$. Then
\Eq{Tak1+}{
  S^T(t,x)=\cl\big(K+T_\alpha(t)\|x\|^\alpha S_0\big)
  \qquad(t\in\R,\,x\in D),
}
where $T_\alpha:\R\to\R$ is the $\alpha$-order Takagi type function defined by
\Eq{Tak2+}{
  T_\alpha(t):=\sum_{n=0}^{\infty}2^{\alpha-n}(d_{\Z}(2^nt))^\alpha
  \qquad(t\in\R).
}}

\begin{proof} To obtain the statement, we can apply \prp{Tak} with the function $\varphi$
defined by $\varphi(x):=\|x\|^\alpha$. Observe that
\Eq{*}{
  \varphi^T(t,x)=\sum_{n=0}^{\infty}\frac{1}{2^n}\varphi\big(2d_{\Z}(2^nt)x\big)
  =\sum_{n=0}^{\infty}2^{\alpha-n} \big(d_{\Z}(2^nt)\big)^\alpha\|x\|^\alpha
  = T_\alpha(t)\|x\|^\alpha
  \qquad(t\in\R,\,x\in D).
}
Therefore, \eq{Tak1+} is a consequence of \eq{Tak1} of \prp{Tak}. 
\end{proof}

\Rem{Tak}{An important particular case is when $\alpha=1$, then $T_1=2T$, where $T$ is the Takagi function 
defined by \eq{Tak} in the introduction. In the case $\alpha=2$ an interesting argument results in a closed form
for $T_2$. Observe that $T_\alpha$ (for any $\alpha>0$) satisfies the functional equation
\Eq{TT}{
  T_\alpha(t)=2^{\alpha}\big(d_{\Z}(t)\big)^\alpha+\frac12 T_\alpha(2t) \qquad(t\in\R).
}
By Banach's fixed-point theorem, this functional equation has a unique solution in the Banach 
space of bounded real functions over the real line (which is equipped with the supremum norm).
Thus $T_\alpha$ is a unique solution to \eq{TT}. On the other hand, for $\alpha=2$, one can easily check
that the 1-periodic function $T_2^*$ defined on $[0,1]$ by $T_2^*(t)=4t(1-t)$ is also a solution of
\eq{TT}, thus we must have $T_2(t)=4t(1-t)$ for $t\in [0,1]$. For further details, see \cite{MakPal13b}.}

\Cor{Tk}{Let $X$ be a normed space, $D\subseteq X$ be a star-shaped set, $S_0\subseteq Y$ be a convex set 
containing $0\in Y$, $K\subseteq Y$ be a convex cone. Define $S:D\to\P_0(Y)$ by $S(x):=K+S_0$. Then
\Eq{Tk}{
  S^T(t,x)=\cl\big(K+2S_0\big)  \qquad(t\in\R,\,x\in D).
}}

\begin{proof} We apply \prp{Tak} for the constant function $\varphi\equiv 1$.
Then \eq{Tak2} yields that $\varphi^T\equiv 2$, whence \eq{Tak1} implies the statement.
\end{proof}

\section{Main Results}
\setcounter{theorem}{0}

The main results of this paper are contained in the following two theorems.
Throughout this section, we assume that $X$ and $Y$ are Hausdorff topological linear spaces.

\Thm{Convex}{Let $D\subseteq X$ be a nonempty convex set and $A,B:(D-D)\to\P_0(Y)$ such that
$0\in A(x)\cap B(x)$ for all $x\in (D-D)$. Denote $\overline{\rec}(B)$, the closure of the
recession cone of $B$ by $K$. Let $F:D\to\P_0(Y)$ be a set-valued mapping which satisfies the 
Jensen-convexity-type inclusion
\Eq{JCV}{
\dfrac{F(x) + F(y)}{2} + A(x-y) \subseteq \cl\bigg(F\bigg(\dfrac{x+y}{2}\bigg) + B(x-y)\bigg) 
   \qquad (x,y\in D).
}
Assume, in addition that $F$ has the following two $K$-boundedness properties.
\begin{enumerate}[(i)]
 \item $F$ is pointwise closedly $K$-lower bounded, i.e., for each $x\in D$, there exists a bounded set 
$H\subseteq Y$ such that $F(x)\subseteq \cl(H+K)$;
 \item $F$ is locally closedly weakly $K$-upper bounded on $D$, i.e., for all $x\in D$, there exist
an open set $U$ containing $x$ and a bounded set $H\subseteq Y$ such that 
$0\in \cl(F(u)+H+K)$ holds for all $u\in U$.
\end{enumerate}
Then $F$ satisfies the convexity type inclusion
\Eq{CV}{
 tF(x)+(1-t)F(y)+A^T(t,x-y)\subseteq \cl\big(F(tx+(1-t)y)+B^T(t,x-y)\big)
   \quad(x,y\in D,\,t\in[0,1]).
}}

\Thm{Concave}{Let $D\subseteq X$ be a nonempty convex set and $A,B:(D-D)\to\P_0(Y)$ such that
$0\in A(x)\cap B(x)$ for all $x\in (D-D)$. Denote $\overline{\rec}(B)$, the closure of the
recession cone of $B$ by $K$. Let $F:D\to\P_0(Y)$ be a set-valued mapping which satisfies the
following Jensen-concavity-type inclusion
\Eq{JCC}{
F\bigg(\dfrac{x+y}{2}\bigg) + A(x-y) \subseteq \cl\bigg(\dfrac{F(x) + F(y)}{2} + B(x-y)\bigg) \qquad (x,y\in D).
}
Assume, in addition that $F$ has the following $K$-convexity and $K$-boundedness properties.
\begin{enumerate}[(i)]
 \item $F$ is pointwise closedly $K$-convex, i.e., $tF(x)+(1-t)F(x)\subseteq \cl(F(x) + K)$ holds 
for each $x\in D$ and for all $t\in[0,1]$;
 \item $F$ is locally closedly $K$-lower bounded, i.e., for each $x\in D$, there exist a neighborhood
$U$ of $x$ and a bounded set $H\subseteq Y$ such that $F(u)\subseteq \cl(H+K)$ for all $u\in D\cap U$.
\end{enumerate}
Then $F$ satisfies the concavity type inclusion
\Eq{CC}{
 F(tx+(1-t)y)+A^T(t,x-y)\subseteq \cl\big(tF(x)+(1-t)F(y)&+B^T(t,x-y)\big)
  \quad(x,y\in D,\,t\in[0,1]).
}}

\Rem{R}{In each of the above theorems the closure operation can be removed from the right hand sides
of the inclusions \eq{JCV}, \eq{CV}, \eq{JCC}, and \eq{CC} if the values of the set-valued map $F$ are 
compact and $B$ has closed values. The closure operation can also be removed from the right hand sides
of \eq{JCV} and \eq{CV} if $F$ has closed values and $B$ is compact valued. This observation also applies
to the corollaries below. Another thing is which is worth mentioning is that if 
$A(0)\subseteq\overline{\rec}(A)$ and $B(0)\subseteq \overline{\rec}(B)$, then, in view of \lem{TT},
the inclusions \eq{CC} and \eq{CV} reduce to \eq{JCC} and \eq{JCV} for the substitution $t=\frac12$,
respectively. Therefore, in this case, under the boundedness and convexity assumptions on $F$, 
\eq{CC} and \eq{CV} are equivalent to \eq{JCC} and \eq{JCV}, respectively. The problem whether 
inclusions \eq{CC} and \eq{CV} are the sharpest possible is an open problem. Results where the
exactness of such estimates were obtained are due to Boros \cite{Bor08}, Makó and Páles \cite{MakPal10b,MakPal13b}.}

The proofs of the above two theorems are described in the next section. In what follows, taking particular 
choices of the set-valued maps $A$, $B$ and using \prp{Tak}, we establish some of the important direct 
consequences of these theorems. They will illustrate how the results recalled in the introduction
are related to our main theorems. 

In the next four corollaries we suppose that $D\subseteq X$ is a nonempty convex set, $K\subseteq Y$ 
is a nonempty closed convex cone, $S_0\subseteq Y$ is a convex set containing $0$ and 
$\varphi:(D-D)\to\R_+$ is a locally upper bounded nonnegative function. Note that, by the convexity of
$D$, the set $(D-D)$ is starshaped, thus \prp{Tak} can be applied. 

The first two corollaries are about approximately and strongly $K$-Jensen convex set-valued mapping, 
respectively.

\Cor{Convex+1}{Assume that $F:D\to\P_0(Y)$ is a pointwise closedly $K$-lower bounded and locally 
closedly weakly $K$-upper bounded set-valued mapping which satisfies
\Eq{JCV+1}{
\dfrac{F(x) + F(y)}{2} \subseteq \cl\bigg(F\bigg(\dfrac{x+y}{2}\bigg) 
    + K + \varphi(x-y)S_0 \bigg) \qquad (x,y\in D).
}
Then 
\Eq{CV+1}{
 tF(x)+(1-t)F(y) \subseteq \cl\big(F(tx+(1-t)y)+K+\varphi^T(t,x-y)S_0\big)\qquad(x,y\in D,\,t\in[0,1]).
}}

\Cor{Convex+2}{Assume that $F:D\to\P_0(Y)$ is a pointwise closedly $K$-lower bounded and locally 
closedly weakly $K$-upper bounded set-valued mapping which satisfies
\Eq{JCV+2}{
\dfrac{F(x) + F(y)}{2} + \varphi(x-y)S_0 \subseteq \cl\bigg(F\bigg(\dfrac{x+y}{2}\bigg) 
  + K \bigg) \qquad (x,y\in D).
}
Then 
\Eq{CV+2}{
 tF(x)+(1-t)F(y) + \varphi^T(t,x-y)S_0 \subseteq \cl\big(F(tx+(1-t)y) + K \big)\qquad(x,y\in D,\,t\in[0,1]).
}}

The next two corollaries are about approximately and strongly $K$-Jensen concave set-valued mapping, 
respectively.

\Cor{Concave+1}{Assume that $F:D\to\P_0(Y)$ is a pointwise closedly $K$-convex and locally closedly 
$K$-lower bounded set-valued mapping which satisfies
\Eq{JCC+1}{
F\bigg(\dfrac{x+y}{2}\bigg) \subseteq \cl\bigg(\dfrac{F(x) + F(y)}{2} 
  + K + \varphi(x-y)S_0 \bigg) \qquad (x,y\in D).
}
Then 
\Eq{CC+1}{
 F(tx+(1-t)y) \subseteq \cl\big(tF(x)+(1-t)F(y) + K + \varphi^T(t,x-y)S_0 \big)\qquad(x,y\in D,\,t\in[0,1]).
}}

\Cor{Concave+2}{Assume that $F:D\to\P_0(Y)$ is a pointwise closedly $K$-convex and locally closedly 
$K$-lower bounded set-valued mapping which satisfies
\Eq{JCC+2}{
F\bigg(\dfrac{x+y}{2}\bigg) + \varphi(x-y)S_0 \subseteq \cl\bigg(\dfrac{F(x) + F(y)}{2} 
    + K\bigg) \qquad (x,y\in D).
}
Then 
\Eq{CC+2}{
 F(tx+(1-t)y)+\varphi^T(t,x-y)S_0 \subseteq \cl\big(tF(x)+(1-t)F(y) + K \big)\qquad(x,y\in D,\,t\in[0,1]).
}}

\begin{proof}[Proof of the Corollaries \ref{CConvex+1}--\ref{CConcave+2}]
Using \thm{Convex} with the set-valued maps $A(u)=0$ and $B(u)=K+\varphi(u)S_0$ (resp., 
$A(u)=\varphi(u)S_0$ and $B(u)=K$) and applying the \prp{Tak}, we obtain \cor{Convex+1} 
(resp., \cor{Convex+2}). Observe that, in both settings, we have that $K\subseteq\overline{\rec}(B)$, 
thus the pointwise closed $K$-lower boundedness and local closed weak $K$-upper boundedness of $F$ imply
its pointwise closed $\overline{\rec}(B)$-lower boundedness and local closed weak $\overline{\rec}(B)$-upper 
boundedness.

Analogously, using \thm{Concave} with the set-valued maps $A(u)=0$ and $B(u)=K+\varphi(u)S_0$ (resp., 
$A(u)=\varphi(u)S_0$ and $B(u)=K$) and applying the \prp{Tak}, we obtain \cor{Concave+1} 
(resp., \cor{Concave+2}). In both settings, we have that $K\subseteq\overline{\rec}(B)$, thus
the pointwise closed $K$-convexity and local closed $K$-lower boundedness of $F$ imply its pointwise 
closed $\overline{\rec}(B)$-convexity and local closed $\overline{\rec}(B)$-lower boundedness.
\end{proof}

\Rem{CC}{The results mentioned and recalled in the introduction can be derived as obvious consequences of
the above corollaries. In the real valued setting, the Bernstein--Doetsch Theorem \cite{BerDoe15}, the 
results of Ng--Nikodem \cite{NgNik93} and Házy--Páles \cite{HazPal04} follow if, in \cor{Convex+1}, we take 
$Y:=\R$, $K:=\R_+$, $S_0:=[-1,0]$, $F(x):=\{f(x)\}$, and $\varphi(x):=0$, $\varphi(x):=\varepsilon$, $\varphi(x):=\varepsilon\|x\|$, respectively. Observe that, in these cases, \prp{Tak} yields 
$\varphi^T(t,x):=0$, $\varphi^T(t,x):=2\varepsilon$, $\varphi^T(t,x):=2\varepsilon T(t)\|x\|$, respectively.
The results of Averna, Cardinali, Nikodem, and Papalini \cite{AveCar90,CarNikPap93,Nik86,Nik87a,Nik87c,Nik89,Pap90} 
and by Borwein \cite{Bor77} that are related to $K$-Jensen convex/concave vector-valued and set-valued mappings
can also be obtained directly. Numerous results obtained for approximate midconvexity by Makó and Páles \cite{MakPal10b,MakPal13b} and by Mureńko, Ja. Tabor, Jó. Tabor, and Żoldak
\cite{MurTabTab12,TabTab09b,TabTab09a,TabTabZol10b,TabTabZol10a} are generalized by Corollaries \ref{CConvex+1}--\ref{CConcave+2} to the vector-valued and set-valued setting. Similarly, using the
explicit form of the function $T_2$ described in \rem{Tak}, one can easily derive the results of Azócar, 
Gimenez, Nikodem and Sanchez \cite{AzoGimNikSan11} and Leiva, Merentes, Nikodem, and Sanchez 
\cite{LeiMerNikSan13} that are related to strongly $K$-Jensen convex real valued and set-valued functions 
from \cor{Convex+2}.}

\section{Proofs of the Two Main Theorems}

\begin{proof}[Proof of \thm{Convex}]
As the first step of the proof of \eq{CV}, we are going to show that, for all $x,y\in D$
there exists a bounded set $H\subseteq Y$ such that, for all $n\geq 0$, $t\in[0,1]$, 
\Eq{CVn}{
 tF(x)+(1-t)F(y) &+ \sum_{k=0}^{n-1}{\dfrac{1}{2^k}A\big(2 d_{\Z}(2^k t)(x-y)\big)} \\
 &\subseteq \cl\bigg(F(tx+(1-t)y) + \dfrac{1}{2^n} H + K +
  \sum_{k=0}^{n-1}{\dfrac{1}{2^k}B\big(2 d_{\Z}(2^k t)(x-y)\big)}\bigg). 
}

Fix $x,y\in D$ arbitrarily. To verify that \eq{CVn} holds, we will proceed by induction on $n.$ 
For the case $n=0$, we have to prove that there exists a bounded set $H$ such that, for all
$t\in[0,1]$, 
\Eq{J0}{
tF(x)+(1-t)F(y)\subseteq \cl\big(F(tx+(1-t)y) + H + K\big).
}
Let $U\in\U(Y)$ and choose a balanced $V\in\U(Y)$ such that $V+V+V\subseteq U$. 
Because $F$ is pointwise closedly $K$-lower bounded, there exist bounded sets $H_x,H_y\subseteq Y$
such that 
\Eq{*}{
F(x) \subseteq \cl(H_x + K)\subseteq V+H_x+K \qquad\mbox{and}\qquad 
F(y) \subseteq \cl(H_y + K) \subseteq V+H_y + K.
} 
Multiplying these inclusions by $t$ and $1-t$, respectively, adding them up side by side, 
and using the convexity of $K$, we obtain
\Eq{J1}{
tF(x)+(1-t)F(y)&\subseteq tV + tH_x + tK + (1-t)V + (1-t)H_y +(1-t)K \\
               &\subseteq V + V + tH_x + (1-t)H_y + K.
}
One can prove that the sets $H_1:=\bigcup_{t\in[0,1]} tH_x$ and $H_2:=\bigcup_{t\in[0,1]} (1-t)H_y$ are
bounded. Thus, inclusion \eq{J1} yields that, for all $t\in[0,1]$,
\Eq{J2}{
tF(x)+(1-t)F(y)\subseteq V + V + H_1 + H_2 + K.
}
On the other hand, applying \lem{wubd}, by the local closed weak $K$-upper boundedness and the compactness
of the segment $[x,y]$, there exist a bounded set $H_0$ such that, for all $t\in[0,1]$,
\Eq{J3}{
 0\in \cl(F(tx+(1-t)y) + H_0 + K)\subseteq V + F(tx+(1-t)y) + H_0 + K.
} 
Now adding up the inclusions \eq{J2} and \eq{J3} side by side, for all $t\in[0,1]$, we obtain
\Eq{*}{
tF(x)+(1-t)F(y)&\subseteq V + V + V+ F(tx+(1-t)y) + H_0 + H_1 + H_2 + K \\
               &\subseteq U + F(tx+(1-t)y) + H_0 + H_1 + H_2 + K.
} 
Therefore,
\Eq{*}{
tF(x)+(1-t)F(y)&\subseteq \bigcap_{U\in\U}\big(U + F(tx+(1-t)y) + H_0 + H_1 + H_2 + K\big)\\
               &=\cl\big(F(tx+(1-t)y) + H_0 + H_1 + H_2 + K\big).
}
Thus, inclusion \eq{J0} follows with $H:=H_0+H_1+H_2$.

Now, suppose that the inclusion \eq{CVn} holds for $n$ and let us prove that it is also valid for
$n+1$. Assume that $t\in\big[0,\frac12\big]$ (the case when $t\in\big[\frac12,1\big]$ is
completely analogous). Then $d_{\Z}(t)=t$ and we can write the left hand side of the inclusion as
\Eq{e11}{
 tF(x)&+(1-t)F(y) + \sum_{k=0}^{n}{\dfrac{1}{2^k}A\big(2 d_{\Z}(2^k t)(x-y)\big)} \\
 &= tF(x)+(1-t)F(y)+ A\big(2t(x-y)\big) 
   + \sum_{k=1}^{n}{\dfrac{1}{2^k}A\big(2 d_{\Z}(2^kt)(x-y)\big)}.
}
We have that 
\Eq{CVC1}{
 (1-t)F(y) \subseteq \frac{1-2t}2 F(y) + \frac12 F(y),
}
and therefore 
\Eq{II1}{
 & tF(x)+(1-t)F(y)+ A\big(2t(x-y)\big) 
   + \sum_{k=1}^{n}{\dfrac{1}{2^k}A\big(2 d_{\Z}(2^kt)(x-y)\big)} \\
 & \subseteq \dfrac{1}{2}\bigg( 2tF(x)+(1-2t)F(y)
   + \sum_{k=0}^{n-1}{\dfrac{1}{2^k}A\big(2 d_{\Z}(2^k(2t))(x-y)\big)} \bigg) 
   + \frac{1}{2}F(y) \!+\! A\big(2d_{\Z}(t)(x-y)\big).
}
Using our inductive hypothesis with $2t$ instead of $t,$ it follows that 
\Eq{II2}{
2tF(x)&+(1-2t)F(y)+ \sum_{k=0}^{n-1} \dfrac{1}{2^k}A\big(2 d_{\Z}(2^k(2t))(x-y)\big) \\
   & \subseteq \cl\bigg(F(2tx + (1-2t)y)+ \frac{1}{2^n}H + K 
   + \sum_{k=0}^{n-1} \dfrac{1}{2^k}B\big(2 d_{\Z}(2^k(2t))(x-y)\big)\bigg).
}
Combining the inclusions \eq{e11}, \eq{II1} and \eq{II2}, we arrive at
\Eq{*}{
 & tF(x)+(1-t)F(y) + \sum_{k=0}^{n}{\dfrac{1}{2^k}A\big(2 d_{\Z}(2^k t)(x-y)\big)} \\
 & \subseteq \dfrac{1}{2}\cl\bigg(F(2tx + (1-2t)y)+ \frac{1}{2^n}H + K 
   + \sum_{k=0}^{n-1} \dfrac{1}{2^k}B\big(2 d_{\Z}(2^k(2t))(x-y)\big)\bigg)
   + \frac{1}{2}F(y) \!+\! A\big(2d_{\Z}(t)(x-y)\big)\\
 & \subseteq \cl\bigg(\dfrac{F(2tx + (1-2t)y)+ F(y)}{2}+ \frac{1}{2^{n+1}}H + K 
     + \sum_{k=0}^{n-1}\dfrac{1}{2^{k+1}}B\big(2d_{\Z}(2^k(2t))(x-y)\big)
     + A\big(2t(x-y)\big)\bigg).
}
By the Jensen-convexity property \eq{JCV} of $F$, we have that
\Eq{*}{
\frac{F(2tx + (1-2t)y) + F(y)}{2} + A\big(2t(x-y)\big)
  \subseteq \cl\bigg(F\bigg(\frac{2tx+(1-2t)y+y}{2}\bigg) + B\big(2t(x-y)\big)\bigg).
}
This and the previous inclusion imply
\Eq{*}{
&tF(x)+(1-t)F(y) + \sum_{k=0}^{n}\dfrac{1}{2^k}A\big(2 d_{\Z}(2^kt)(x-y)\big) \\
& \subseteq \cl\bigg(\!\cl\bigg(F\bigg(\frac{2tx+(1-2t)y+y}{2}\bigg) 
     + B\big(2t(x-y)\big)\bigg) + \frac{1}{2^{n+1}}H + K 
     + \sum_{k=0}^{n-1}\dfrac{1}{2^{k+1}}B\big(2d_{\Z}(2^k(2t))(x-y)\big)\bigg)\\
& = \cl\bigg(F\big(tx+(1-t)y\big) + \frac{1}{2^{n+1}}H + K 
  + \sum_{k=0}^{n}\dfrac{1}{2^{k}}B\big(2d_{\Z}(2^k t)(x-y)\big)\bigg).
}
Now, we can conclude that inclusion \eq{CVn} holds for all $n\geq 0.$

To complete the proof of the theorem, let $t\in[0,1]$ be also fixed and apply \lem{Inc} to the
sequences of sets and numbers defined for $n\geq0$ as
\Eq{*}{
 A_n &:= tF(x)+(1-t)F(y) + \sum_{k=0}^{n-1}\dfrac{1}{2^k}A(2 d_{\Z}(2^k t)(x-y)), \\
 B_n &:= F\big(tx+(1-t)y\big) + \sum_{k=0}^{n-1}\dfrac{1}{2^{k}}B\big(2d_{\Z}(2^k t)(x-y)\big), \\
 \varepsilon_n &:= \frac{1}{2^{n}}.
}
Then, with these notations, inclusion \eq{CVn} is equivalent to \eq{Incn}.
On the other hand, by the assumption that $0\in A(u)\cap B(u)$ for all $u\in (D-D)$,
it easily follows that $(A_n)$ and $(B_n)$ are nondecreasing sequences of subsets of $Y$.

We will show that $K\subseteq\bigcap_{n=0}^\infty\overline{\rec}(B_n)$. For this, it suffices to prove, that
for all neighborhood $U$ of zero and for all $n\in\N$, the inclusion $K\subseteq U+\rec(B_n)$ is valid.
Let $n\in\N$ and let $U$ be an arbitrary neighborhood of zero. Then choose a neighborhood $V$ of zero such that
$\sum_{k=0}^{n-1} 2^{-k} V\subseteq U$. Then, $K\subseteq V+\rec(B(u))$ for all $u\in D-D$. Hence, using the properties of the recession cones established in \lem{rec}, we obtain
\Eq{*}{
  K=\sum_{k=0}^{n-1} \frac1{2^k}K
   &\subseteq\sum_{k=0}^{n-1} \frac1{2^k}\Big(V+\rec\big(B\big(2d_{\Z}(2^k t)(x-y)\big)\big)\Big) \\
   &\subseteq U+\rec\bigg(\sum_{k=0}^{n-1}\frac1{2^k}\big(B\big(2d_{\Z}(2^k t)(x-y)\big)\bigg) \\
   &\subseteq U+\rec\bigg(F\big(tx+(1-t)y\big) 
     +\sum_{k=0}^{n-1}\frac1{2^k}\big(B\big(2d_{\Z}(2^k t)(x-y)\big)\bigg),
}
which proves the inclusion $K\subseteq U+\rec(B_n)$. Therefore $K\subseteq\bigcap_{n=0}^\infty\overline{\rec}(B_n)$.

Now we are in the position to apply \lem{Inc}. Hence \eq{Inc} holds, in other words, we obtain that
\Eq{*}{
  \cl\bigg(\bigcup_{n=1}^\infty \bigg(tF(x)+(1-t)F(y)
     &+ \sum_{k=0}^{n-1}\dfrac{1}{2^k}A(2 d_{\Z}(2^k t)(x-y))\bigg)\bigg)\\
  &\subseteq \cl\bigg(\bigcup_{n=1}^\infty \bigg(F\big(tx+(1-t)y\big)
     + \sum_{k=0}^{n-1}\dfrac{1}{2^{k}}B\big(2d_{\Z}(2^k t)(x-y)\big)\bigg)\bigg).
}
Now applying the equality in \eq{cl+} to evaluate the left and right hand sides of the above inclusion, 
it follows that
\Eq{*}{
  \cl\bigg(tF(x)+(1-t)F(y) &+ \cl\bigg(\bigcup_{n=1}^\infty
     \sum_{k=0}^{n-1}\dfrac{1}{2^k}A(2 d_{\Z}(2^k t)(x-y))\bigg)\bigg)\\
  &\subseteq \cl\bigg(F\big(tx+(1-t)y\big) +\cl\bigg(\bigcup_{n=1}^\infty 
     \sum_{k=0}^{n-1}\dfrac{1}{2^{k}}B\big(2d_{\Z}(2^k t)(x-y)\big)\bigg)\bigg),
}
which is equivalent to the inclusion \eq{CV} to be proved.
\end{proof}

\begin{proof}[Proof of \thm{Concave}]
To prove \eq{CC}, we are going to show first that, for all $x,y\in D$, there exists a bounded set
$H\subseteq Y$ such that, for all $n\geq 0$ and $t\in[0,1]$, 
\Eq{CCn}{
 F(tx+(1-t)y) &+ \sum_{k=0}^{n-1}\dfrac{1}{2^k}A\big(2 d_{\Z}(2^k t)(x-y)\big) \\[-2mm]
 &\subseteq \cl\bigg(tF(x)+(1-t)F(y) + \dfrac{1}{2^n} H + K +
  \sum_{k=0}^{n-1}\dfrac{1}{2^k}B\big(2 d_{\Z}(2^k t)(x-y)\big)\bigg).
}

Let $x,y\in D$ be fixed. To verify that \eq{CCn} holds, we will proceed by induction over $n.$ 

The bounded set $H$ will be constructed so that \eq{CCn} be valid for $n=0$, that is, for all
$t\in[0,1]$, the following condition holds:
\Eq{CC0}{
  F(tx+(1-t)y)\subseteq \cl\big(tF(x)+ (1-t)F(y) + H + K\big).
}
In view of \lem{lbd}, the local closed $K$-lower boundedness of $F$ and the compactness of the 
segment $[x,y]$ imply that there exists a bounded set $H_0\subseteq Y$ such that 
\Eq{H0}{
F(tx+(1-t)y)\subseteq \cl(H_0 + K)   \qquad  (t \in[0,1]).
}
On the other hand, the sets $F(x)$ and $F(y)$ being nonempty, we can choose two elements
$u\in F(x)$ and $v\in F(y)$. Then
\Eq{Hxy}{
0 \in F(x)-u \qquad\mbox{and}\qquad 0 \in F(y)-v.
}    
Multiplying the two inclusions in \eq{Hxy} by $t$ and $(1-t)$, respectively, and adding them up
to the inclusion \eq{H0}, for $t\in[0,1]$, we obtain that 
\Eq{*}{
F(tx+(1-t)y) &\subseteq tF(x)+(1-t)F(y) - tu - (1-t)v + \cl(H_0 + K)\\
            &\subseteq \cl\big(tF(x)+(1-t)F(y) -[u,v] + H_0 + K\big).
}
Therefore, \eq{CC0} holds with $H := H_0-[u,v],$ which is obviously bounded.

Now, suppose that \eq{CCn} is valid for $n$ and let us prove its validity for $n+1.$
Assume that $t\in\big[0,\frac12\big]$. Observe that then $d_{\Z}(t)=t$. 
Let us start evaluating the left side of the inclusion to be proved. 
\Eq{I1}{
 F(tx&+(1-t)y) + \sum_{k=0}^{n}{\dfrac{1}{2^k}A\big(2 d_{\Z}(2^k t)(x-y)\big)} \\
 & = F(tx+(1-t)y) + A(2d_{\Z}(t)(x-y)) + \sum_{k=1}^{n}{\dfrac{1}{2^k}A\big(2 d_{\Z}(2^k t)(x-y)\big)} \\
 & = F\bigg(\dfrac{2tx+(1-2t)y + y}{2}\bigg) + A\big(2t(x-y)\big)
   + \dfrac{1}{2}\sum_{k=0}^{n-1}{\dfrac{1}{2^k}A\big(2 d_{\Z}(2^k (2t))(x-y)\big)}.
}
By the Jensen-concavity property \eq{JCC} of $F$, we get the inclusion
\Eq{I2}{
 F\bigg(\dfrac{2tx+(1-2t)y + y}{2}\bigg) + A\big(2t(x-y)\big) 
  \subseteq \cl\bigg(\dfrac{F(2tx+(1-2t)y) + F(y)}{2} + B\big(2t(x-y)\big)\bigg).
}
Combining \eq{I1} and \eq{I2}, we obtain
\Eq{I3}{
 F&(tx+(1-t)y) + \sum_{k=0}^{n}\dfrac{1}{2^k}A\big(2 d_{\Z}(2^k t)(x-y)\big) \\
  &\subseteq \cl\bigg(\dfrac{F(2tx+(1-2t)y) + F(y)}{2} + B\big(2t(x-y)\big)\bigg)
   +\dfrac{1}{2}\sum_{k=0}^{n-1}\dfrac{1}{2^k}A\big(2 d_{\Z}(2^k (2t))(x-y)\big) \\
  & \subseteq \cl\bigg(\dfrac12\bigg(F(2tx+(1-2t)y)
   +\sum_{k=0}^{n-1}\dfrac{1}{2^k}A\big(2 d_{\Z}(2^k (2t))(x-y)\big)\bigg)
   + \dfrac12 F(y)+B\big(2d_{\Z}(t)(x-y)\big)\bigg).
}
Now, using our inductive hypothesis with $2t$ instead of $t$, we get 
\Eq{I4}{  
 F(2tx&+(1-2t)y)+\sum_{k=0}^{n-1}\dfrac{1}{2^k}A\big(2 d_{\Z}(2^k (2t))(x-y)\big) \\ 
&\subseteq \cl\bigg(2tF(x)+(1-2t)F(y) + \dfrac{1}{2^{n}}H + K +
\sum_{k=0}^{n-1}\dfrac{1}{2^k}B\big(2 d_{\Z}(2^k (2t))(x-y)\big)\bigg).
}
Inserting \eq{I4} into \eq{I3} and using that $(1-2t)F(y)+F(y)\subseteq \cl\big((2-2t)F(y)+K\big)$
(which is a consequence of the pointwise closed $K$-convexity of $F$), it follows that
\Eq{*}{
 F(tx&+(1-t)y) + \sum_{k=0}^{n}{\dfrac{1}{2^k}A\big(2 d_{\Z}(2^k t)(x-y)\big)} \\
& \subseteq \cl\bigg(\dfrac12\cl\bigg(2tF(x)+(1-2t)F(y) + \dfrac{1}{2^{n}}H + K +
\sum_{k=0}^{n-1}\dfrac{1}{2^k}B\big(2 d_{\Z}(2^k (2t))(x-y)\big)\bigg)\\
   &\qquad + \dfrac12 F(y)+B\big(2d_{\Z}(t)(x-y)\big)\bigg)\\
&\subseteq \cl\bigg(\dfrac{1}{2}\big(2t F(x) + (1-2t) F(y) + F(y) + K\big) + \dfrac{1}{2^{n+1}}H 
 +\sum_{k=0}^{n}\dfrac{1}{2^{k}}B\big(2 d_{\Z}(2^k t)(x-y)\big)\bigg) \\
&\subseteq \cl\bigg(tF(x)+(1-t)F(y)+ \dfrac{1}{2^{n+1}}H + K +
\sum_{k=0}^{n}\dfrac{1}{2^{k}}B\big(2 d_{\Z}(2^k t)(x-y)\big)\bigg). 
}
This completes the proof of the induction and hence \eq{CCn} holds for all $n\geq0$.

Now we are going to use \lem{Inc}, so that, for a fixed $t\in [0,1]$, we define the sequences 
\Eq{*}{
A_n &:= F(tx+(1-t)y) + \sum_{k=0}^{n-1}{\dfrac{1}{2^k}A\big(2 d_{\Z}(2^k t)(x-y)\big)}, \\
B_n &:= tF(x)+(1-t)F(y) + \sum_{k=0}^{n-1}{\dfrac{1}{2^k}B\big(2 d_{\Z}(2^k t)(x-y)\big)}, \\
\varepsilon_n &= \frac1{2^n}.
}

Then, inclusion \eq{CCn}, with sequences $(A_n),(B_n)$ and $(\varepsilon_n)$ defined 
above is equivalent to \eq{Incn}. One can see that the sequences $(A_n)$ and $(B_n)$ 
are nondecreasing and $K\subseteq\bigcap_{n=0}^\infty\overline{\rec}(B_n)$ also holds
(due to a similar argument that was followed in the proof of \thm{Convex}).
Thus, by the \lem{Inc}, it follows that
\Eq{*}{
  \cl\bigg(\bigcup_{n=0}^{\infty}F(tx+(1-t)y) 
   &+ \sum_{k=0}^{n-1}{\dfrac{1}{2^k}A\big(2 d_{\Z}(2^kt)(x-y)\big)}\bigg)\\
  &\subseteq \cl\bigg(\bigcup_{n=0}^{\infty}tF(x)+(1-t)F(y) 
   + \sum_{k=0}^{n-1}{\dfrac{1}{2^k}B\big(2 d_{\Z}(2^k t)(x-y)\big)}\bigg).
}
Now, similarly as in the proof of \thm{Convex}, using \eq{cl+}, this relation implies the
desired inclusion \eq{CC}.
\end{proof}

%\nocite{Kuc85}
%\bibliography{publ,funcequ,control}

\begin{thebibliography}{10}

\bibitem{AveCar90}
A.~Averna and T.~Cardinali, \emph{On the concepts of {$K$}-convexity
  [{$K$}-concavity] and {$K$}-convexity{$^*$} [{$K$}-concavity{$^*$}]}, Riv.
  Mat. Univ. Parma (4) \textbf{16} (1990), no.~1-2, 311–330. \MR{1105752
  (92h:26031)}

\bibitem{AzoGimNikSan11}
A.~Azócar, J.~Giménez, K.~Nikodem, and J.~L. Sánchez, \emph{On strongly
  midconvex functions}, Opuscula Math. \textbf{31} (2011), no.~1, 15–26.
  \MR{2739838 (2011k:26009)}

\bibitem{BerDoe15}
F.~Bernstein and G.~Doetsch, \emph{Zur {T}heorie der konvexen {F}unktionen},
  Math. Ann. \textbf{76} (1915), no.~4, 514–526. \MR{1511840}

\bibitem{Bor08}
Z.~Boros, \emph{{An inequality for the {T}akagi function}}, Math. Inequal.
  Appl. \textbf{11} (2008), no.~4, 757–765. \MR{2009f:39047}

\bibitem{Bor77}
J.M. Borwein, \emph{Multivalued convexity and optimization: a unified approach
  to inequality and equality constraints}, Math. Programming \textbf{13}
  (1977), no.~2, 183–199. \MR{0451166 (56 \#9453)}

\bibitem{CarNikPap93}
T.~Cardinali, K.~Nikodem, and F.~Papalini, \emph{Some results on stability and
  on characterization of {$K$}-convexity of set-valued functions}, Ann. Polon.
  Math. \textbf{58} (1993), no.~2, 185–192. \MR{1239022 (94g:26022)}

\bibitem{Haz05a}
A.~Házy, \emph{On approximate {$t$}-convexity}, Math. Inequal. Appl.
  \textbf{8} (2005), no.~3, 389–402. \MR{2148233 (2006c:26019)}

\bibitem{Haz07b}
A.~Házy, \emph{On the stability of {$t$}-convex functions}, Aequationes Math.
  \textbf{74} (2007), no.~3, 210–218. \MR{2376448 (2008j:26012)}

\bibitem{HazPal04}
A.~Házy and Zs. Páles, \emph{On approximately midconvex functions}, Bull.
  London Math. Soc. \textbf{36} (2004), no.~3, 339–350. \MR{2038721
  (2004j:26020)}

\bibitem{HazPal05}
A.~Házy and Zs. Páles, \emph{On approximately {$t$}-convex functions}, Publ. Math. Debrecen
  \textbf{66} (2005), no.~3-4, 489–501. \MR{2137784 (2006c:26023)}

\bibitem{HazPal09}
A.~Házy and Zs. Páles, \emph{On a certain stability of the {H}ermite-{H}adamard inequality},
  Proc. R. Soc. Lond. Ser. A Math. Phys. Eng. Sci. \textbf{465} (2009),
  no.~2102, 571–583. \MR{2471774 (2009k:39033)}

\bibitem{Kuc85}
M.~Kuczma, \emph{An {I}ntroduction to the {T}heory of {F}unctional {E}quations
  and {I}nequalities}, Prace Naukowe Uniwersytetu Śląskiego w Katowicach,
  vol. 489, Państwowe Wydawnictwo Naukowe — Uniwersytet Śląski,
  Warszawa–Kraków–Katowice, 1985, 2nd edn. (ed. by A. Gilányi),
  Birkhäuser, Basel, 2009. \MR{0788497 (86i:39008), MR 2467621}

\bibitem{Lac99}
M.~Laczkovich, \emph{The local stability of convexity, affinity and of the
  {J}ensen equation}, Aequationes Math. \textbf{58} (1999), 135–142.
  \MR{1714327 (2001d:39028)}

\bibitem{LeiMerNikSan13}
H.~Leiva, N.~Merentes, K.~Nikodem, and J.~L. Sánchez, \emph{Strongly convex
  set-valued maps}, J. Global Optim. \textbf{57} (2013), 695–705.
  \MR{3119375}

\bibitem{MakPal10b}
J.~Makó and Zs. Páles, \emph{Approximate convexity of {T}akagi type
  functions}, J. Math. Anal. Appl. \textbf{369} (2010), no.~2, 545–554.
  \MR{2651700 (2011k:26011)}

\bibitem{MakPal12b}
J.~Makó and Zs. Páles, \emph{Implications between approximate convexity properties and
  approximate {H}ermite-{H}adamard inequalities}, Cent. Eur. J. Math.
  \textbf{10} (2012), no.~3, 1017–1041. \MR{2902231}

\bibitem{MakPal12c}
J.~Makó and Zs. Páles, \emph{Korovkin type theorems and approximate {H}ermite-{H}adamard
  inequalities}, J. Approx. Theory \textbf{164} (2012), no.~8, 1111–1142.
  \MR{2935448}

\bibitem{MakPal13b}
J.~Makó and Zs. Páles, \emph{On approximately convex {T}akagi type functions}, Proc. Amer.
  Math. Soc. \textbf{141} (2013), no.~6, 2069–2080. \MR{3034432}

\bibitem{MurTabTab12}
A.~Mureńko, Ja. Tabor, and Jó. Tabor, \emph{Applications of de {R}ham
  {T}heorem in approximate midconvexity}, J. Diff. Equat. Appl. \textbf{18}
  (2012), no.~3, 335–344. \MR{2901825}

\bibitem{NgNik93}
C.~T. Ng and K.~Nikodem, \emph{On approximately convex functions}, Proc. Amer.
  Math. Soc. \textbf{118} (1993), no.~1, 103–108. \MR{1159176 (93f:26006)}

\bibitem{Nik86}
K.~Nikodem, \emph{Continuity of {$K$}-convex set-valued functions}, Bull.
  Polish Acad. Sci. Math. \textbf{34} (1986), no.~7-8, 393–400. \MR{874882
  (88a:26040)}

\bibitem{Nik87a}
K.~Nikodem, \emph{On concave and midpoint concave set-valued functions}, Glas.
  Mat. Ser. III \textbf{22(42)} (1987), no.~1, 69–76. \MR{940094 (89g:39017)}

\bibitem{Nik87c}
K.~Nikodem, \emph{On midpoint convex set-valued functions}, Aequationes Math.
  \textbf{33} (1987), no.~1, 46–56. \MR{901799 (88h:90171)}

\bibitem{Nik89}
K.~Nikodem, \emph{{$K$}-convex and {$K$}-concave set-valued functions}, Zeszyty
  Nauk. Politech. Łódz. Mat. (Łódz) \textbf{559} (1989), 1–75, (Rozprawy
  Nauk. 114).

\bibitem{Pap90}
F.~Papalini, \emph{The {$K$}-midpoint {$^*$} convexity [concavity] and lower
  [upper] {$K$}-semicontinuity of a multifunction}, Riv. Mat. Univ. Parma (4)
  \textbf{16} (1990), no.~1-2, 149–159 (1991). \MR{1105736 (92h:26032)}

\bibitem{Pol66}
B.~T. Polyak, \emph{Existence theorems and convergence of minimizing sequences
  for extremal problems with constraints}, Dokl. Akad. Nauk SSSR \textbf{166}
  (1966), 287–290. \MR{33 \#6466}

\bibitem{Rud91}
W.~Rudin, \emph{Functional {A}nalysis}, second ed., International Series in
  Pure and Applied Mathematics, McGraw-Hill Inc., New York, 1991. \MR{1157815
  (92k:46001)}

\bibitem{Rad60}
H.~Rådström, \emph{One-parameter semigroups of subsets of a real linear
  space}, Ark. Mat. \textbf{4} (1960), 87–97. \MR{0146280 (26 \#3802)}

\bibitem{TabTab09b}
Ja. Tabor and Jó. Tabor, \emph{Generalized approximate midconvexity}, Control
  Cybernet. \textbf{38} (2009), no.~3, 655–669. \MR{2650358 (2011f:52002)}

\bibitem{TabTab09a}
Ja. Tabor and Jó. Tabor, \emph{Takagi functions and approximate midconvexity}, J. Math. Anal.
  Appl. \textbf{356} (2009), no.~2, 729–737.

\bibitem{TabTabZol10b}
Ja. Tabor, Jó. Tabor, and M.~Żołdak, \emph{Approximately convex functions on
  topological vector spaces}, Publ. Math. Debrecen \textbf{77} (2010),
  115–123. \MR{2675738 (2011f:26009)}

\bibitem{TabTabZol10a}
Ja. Tabor, Jó. Tabor, and M.~Żołdak, \emph{Optimality estimations for approximately midconvex functions},
  Aequationes Math. \textbf{80} (2010), 227–237. \MR{2736954 (2011j:26020)}

\bibitem{Tru84}
L.~I. Trudzik, \emph{Continuity properties of vector-valued convex functions},
  J. Austral. Math. Soc. Ser. A \textbf{36} (1984), no.~3, 404–415.
  \MR{733912 (85d:46062)}

\end{thebibliography}
%\bibliographystyle{amsplain}

\providecommand{\bysame}{\leavevmode\hbox to3em{\hrulefill}\thinspace}
\providecommand{\MR}{\relax\ifhmode\unskip\space\fi MR }
\providecommand{\MRhref}[2]{%
  \href{http://www.ams.org/mathscinet-getitem?mr=#1}{#2}
}
\providecommand{\href}[2]{#2}

\end{document}